\documentclass{amsart}[11pt]
\usepackage{amssymb, mathrsfs, amsmath}
\newtheorem{Theorem}{Theorem}

\newtheorem{Lemma}{Lemma}

\title{Radon Transform on spheres and Generalized Bessel function associated with dihedral groups} 
\begin{document}
\maketitle
\centerline{N. Demni\footnote{IRMAR, Universit\'e de Rennes 1, Campus de Beaulieu, 35042 Rennes Cedex, France. \\ 
E-mail: nizar.demni@univ-rennes1.fr.  
}}
\begin{abstract}
Motivated by Dunkl operators theory, we consider a generating series involving a modified Bessel function and a Gegenbauer polynomial, that generalizes a known series already considered by L. Gegenbauer. We actually use inversion formulas for Fourier and Radon transforms on spheres to derive a closed formula for this series when the parameter of the Gegenbauer polynomial is a positive integer. As a by-product, we get a relatively simple integral representation for the generalized Bessel function associated with dihedral groups ${\it D}_{n}, n \geq 2$ when both multiplicities sum to an integer. In particular, we recover a previous result obtained for ${\it D}_4$ and we give a special interest to ${\it D}_6$. Finally, we derive similar results for odd dihedral groups.    
\end{abstract}

\section{Motivation}
The dihedral group ${\it D}_n$ of order $n \geq 2$ is defined as the group of regular $n$-gone preserving-symmetries (\cite{Dunkl}). It figures among reflection groups associated with irreducible root systems yet ceases to be crystallographic unless $n=2,3,4,6$. Nevertheless the theory of rational Dunkl operators introduced in the late eighties associates to reduced non necessarily crystallographic root systems generalized Bessel functions that extend spherical functions on symmetric spaces of Euclidean type from a discrete to a continuous range of multiplicities (see Ch.I in \cite{Chy}). In fact, the radial part of the Laplace-Beltrami operator on these symmetric spaces fits the reflection group-invariant part of the Dunkl Laplacian with special multiplicities (this fact holds for radial parts of a more general class of differential operators). However, they are not easy to handle unlike spherical functions, except possibly in lower ranks. In fact, they are expressed for the four infinite series of irreducible root systems as multivariate hypergeometric series defined via Jack polynomials (\cite{Dem03}). Nonetheless, Jack polynomials may be expressed by means of Gegenbauer polynomials in ranks one and two (\cite{Mang}). Moreover, probabilistic considerations led to the following expression for the generalized Bessel function associated with dihedral systems (\cite{Dem01}). Let $n=2p, p \geq 1$ and let $D_k^W$ denote the generalized Bessel function depending in this case on the two real variables $x = \rho e^{i\phi}, y = re^{i\theta}, \rho, r \geq 0, \phi, \theta \in [0,\pi/2p]$. Then
\begin{equation}
\label{GBFE}
D_k^W(\rho,\phi,r,\theta) = c_{p,k}\left(\frac{2}{r\rho}\right)^{\gamma} \sum_{j \geq 0}{\it I}_{2jp +\gamma}(\rho r)p_j^{l_1, l_0}(\cos(2p\phi))p_j^{l_1, l_0}(\cos(2p\theta))
\end{equation}
where 
\begin{itemize}
\item $k=(k_0,k_1)$ is a positive-valued multiplicity function, $l_i = k_i - 1/2, i \in \{1,2\}$, $\gamma = p(k_0+k_1)$. 
\item  ${\it I}_{2jp+\gamma}, p_j^{l_1,l_0}$ are the modified Bessel function of index $2jp+\gamma$ and the $j$-th orthonormal Jacobi polynomial of parameters $l_1,l_0$ respectively (the orthogonality (Beta) measure need not to be normalized here. In fact, the normalization only alters the constant $c_{p,k}$ below). 
\item The constant $c_{p,k}$ depends on $p, k$ and is such that $D_k^W(0,y) = 1$ for all $y = (r,\theta) \in [0,\infty) \times [0,\pi/2p]$ (see \cite{Dem02})
\begin{equation*}
c_{p,k} = 2^{k_0+k_1}\frac{\Gamma(p(k_1+k_0)+1)\Gamma(k_1+1/2)\Gamma(k_0+1/2)}{\Gamma(k_0+k_1+1)}.
\end{equation*}
\item A similar formula holds for odd dihedral systems (see the fourth section). 
\end{itemize} 
Once this relatively simple formula was obtained, the special case $p=2$ was the main object of a subsequent paper (\cite{Dem02}), aiming to work out the series displayed in \eqref{GBFE}. The main achievement was then realized when $\gamma = 2(k_0+k_1)$ is an even integer and according to Corollary 1.2 in \cite{Dem02}
\begin{equation*}
D_k^W(\rho,\phi,r,\theta) = \int \int i_{(\gamma -1)/2} \left(\rho r \sqrt{\frac{1+z_{2\phi,2\theta}(u,v)}{2}}\right)\mu^{l_1}(du)\mu^{l_0}(dv)
\end{equation*}
where 
\begin{equation*}
i_{\alpha}(x) := \sum_{m=0}^{\infty} \frac{1}{(\alpha+1)_m m!} \left(\frac{x}{2}\right)^{2m}, \, \alpha > -1,
\end{equation*}
is the normalized modified Bessel function (\cite{Dunkl}). In this paper, we shall see that this achievement is not specific to the value $p=2$ but rather extends to all $p \geq 1$ provided that $k_0 + k_1$ is a positive integer and is even related to geometrical considerations on spheres that considerably avoid tedious computations performed in \cite{Dem02}. This is seen as follows: start with Dijksma-Koornwinder product formula for Jacobi polynomials (\cite{Dij}) which may be written in the following way (\cite{Dem02}): 
\begin{equation*}
c(\alpha,\beta) p_j^{\alpha,\beta}(\cos2\phi)p_j^{\alpha,\beta}(\cos2\theta) = (2j+\alpha+\beta+1)\int\int C_{2j}^{\alpha+\beta+1}(z_{\phi, \theta}(u,v))\mu^{\alpha}(du)\mu^{\beta}(dv)
\end{equation*}  
where $\alpha, \beta > -1/2$, 
\begin{equation*}
c(\alpha,\beta) = 2^{\alpha+\beta+1}\frac{\Gamma(\alpha+1)\Gamma(\beta+1)}{\Gamma(\alpha+\beta+1)},
\end{equation*}
\begin{equation*}
z_{\phi,\theta}(u,v) = u\cos \theta\cos\phi + v \sin\theta\sin\phi,
\end{equation*}
and $\mu^{\alpha}$ is the symmetric Beta probability measure whose density is given by
\begin{equation*}
\mu^{\alpha}(du) = \frac{\Gamma(\alpha+1)}{\sqrt{\pi}\Gamma(\alpha+1/2)}(1-u^2)^{\alpha-1/2} {\bf 1}_{[-1,1]}(u)du, \quad \alpha > -1/2.
\end{equation*}
Next, invert the order of integration in \eqref{GBFE} to see that
\begin{align}\label{S1}
D_k^W(\rho,\phi,r,\theta) \propto \int \int \left(\frac{2}{r\rho}\right)^{\gamma} & \sum_{j \geq 0}(2j+k_0+k_1){\it I}_{2jp +\gamma}(\rho r)C_{2j}^{k_0+k_1}(z_{p\phi,p\theta}(u,v)) \nonumber
\\&   \mu^{\alpha}(du)\mu^{\beta}(dv)
\end{align}
where the notation $\propto$ means that equality holds up to a constant factor. Now, the integrand in \eqref{S1} is obviously the sum of the following series 
\begin{equation}\label{Series}
f_{\nu,p}^{\pm}(R,\cos \zeta) := \left(\frac{2}{R}\right)^{p\nu}\sum_{j \geq 0}(\pm 1)^j(j+\nu){\it I}_{p(j +\nu)}(R)C_{j}^{\nu}(\cos \zeta) 
\end{equation}
where we set $\nu := k_0+k_1, R := \rho r$ and $\cos \zeta := \cos \zeta(u,v) = z_{p\phi, p\theta}(u,v)$.

Note actually that closed formulas for $f_{\nu,1}^{\pm}$ are due to L. Gegenbauer (equations (4), (5), p.369 in \cite{Watson})
\begin{equation*}
\left(\frac{2}{r\rho}\right)^{\gamma} \sum_{j \geq 0}(\pm 1)^j(j+\gamma){\it I}_{j +\gamma}(\rho r)C_{j}^{\gamma}(\cos \zeta) = \frac{1}{\Gamma(\gamma)}e^{\pm \rho r \cos \zeta}
\end{equation*}
and were used in \cite{Dem02}. Note also that $f_{\nu,p}^-(R,\cos\zeta) = f_{\nu,p}^+(R, -\cos\zeta)$ which follows from $C_j^{\nu}(\cos\zeta) = (-1)^jC_j(-\cos\zeta)$. Our main result is then stated as
\begin{Theorem}
Assume $\nu$ is an integer and $\nu \geq 1$, then 
\begin{equation*}
\left(\frac{R}{2}\right)^{p\nu}f_{\nu,p}^{+}(R,\cos \zeta) = \frac{1}{2^{\nu}(\nu-1)!} \left[-\frac{1}{\sin \zeta} \frac{d}{d\zeta}\right]^{\nu} \frac{1}{p} \sum_{s=1}^p e^{ R\cos[(\zeta+2\pi s)/p]}
\end{equation*}
and 
\begin{equation*}
\left(\frac{R}{2}\right)^{p\nu}f_{\nu,p}^{-}(R,\cos \zeta) = \frac{1}{2^{\nu}(\nu-1)!} \left[\frac{1}{\sin \zeta} \frac{d}{d\zeta}\right]^{\nu} \frac{1}{p} \sum_{s=1}^p e^{ R\cos[(\zeta+2\pi s+\pi)/p]}.
\end{equation*}
\end{Theorem}
We shall write two different proofs of this result. The first one rely on interpreting the sequence
\begin{equation*}
(\pm 1)^j{\it I}_{p(j+\nu)}(R), \, j \geq 0
\end{equation*}
for fixed $R$ as the Gegenbauer-Fourier coefficients of $\zeta \mapsto f_{\nu,p}^{\pm}(R, \cos \zeta)$ corresponding to the Gegenbauer-Fourier transform studied in \cite{Abou}. Thus, deriving closed formulas for $f_{\nu,p}^{\pm}$ when $\nu$ is a positive integer amounts to appropriately use inversion formulas for Fourier and the so-called Radon transforms. The second one do not make use of the Radon transform on spheres and rather uses the formula (\cite{Erd0})
\begin{equation*}
\frac{1}{2^{\nu-1}\Gamma(\nu)} \frac{d^{\nu}}{dz^{\nu}}T_{j+\nu}(z) = (j+\nu)C_j^{\nu}(z),
 \end{equation*}
where $T_j$ is the $j$-th Tchbycheff polynomial of the first kind (\cite{Erd}). The remainder of the paper is devoted to a focus on particular values of $p \geq 2$. When $p=2$, we shall use  
\begin{equation*}
\cos (\zeta/2) = \sqrt{\frac{1+\cos \zeta}{2}}, \quad  \zeta \in [0,\pi].
\end{equation*} 
together with appropriate formulas for modified Bessel functions in order to recover Corollary 1.2. in \cite{Dem02}, while when $p = 3$ we shall solve a special cubic equation when $p=3$ relying on results from analytic function theory rather than Cardan formulas. The required solution is then expressed by means of Gauss hypergeometric functions (\cite{Hille}) yielding therefore a somehow explicit formula for the series \eqref{S1}, though much more complicated than the one derived for $p=2$. More generally, one needs to write down the inverse of the $p$-th Tchebycheff polynomial of the first kind restricted to the interval $[\cos(\pi/p),1]$, which is by the virtue of Galois theory not possible using radicals. The paper is closed with adapting our method to odd dihedral groups, in particular to ${\it D}_3$ thereby exhausting the list of dihedral groups that are Weyl groups,\footnote{The value $p=1$ corresponds to the product group $(\mathbb{Z}_2)^2$ for which $D_k^W$ is namely a product of normalized modified Bessel functions.} and with an explanation of the occurrence of the Radon transform on spheres in our framework. 

\section{Proofs of the main result}
\subsection{First proof via Radon transform} Recall the orthogonality relation for Gegenbauer polynomials (\cite{Dunkl}): 
\begin{align*}
\int_0^{\pi} C_j^{\nu}(\cos \zeta)C_m^{\nu}(\cos \zeta) (\sin \zeta)^{2\nu} d\zeta &= \delta_{jm} \frac{\pi\Gamma(j+2\nu)2^{1-2\nu}}{\Gamma^2(\nu)(j+\nu)j!} 
\\&= \delta_{jm} \frac{\pi 2^{1-2\nu}\Gamma(2\nu)}{(j+\nu) \Gamma^2(\nu)}C_j^{\nu}(1) 
\\& = \delta_{jm} \nu \frac{\sqrt{\pi}\Gamma(\nu+1/2)}{\Gamma(\nu+1)}\frac{C_j^{\nu}(1)}{(j+\nu)}
\end{align*}
where we used Gauss duplication formula (\cite{Dunkl})
\begin{equation*}
\sqrt{\pi}\Gamma(2\nu) = 2^{2\nu-1}\Gamma(\nu)\Gamma(\nu+1/2),
\end{equation*}
and the special value (\cite{Dunkl})
\begin{equation*}
C_j^{\nu}(1) = \frac{(2\nu)_j}{j!} .
\end{equation*}
Equivalently, if $\mu^{\nu}(d\cos \zeta)$ is the image of $\mu^{\nu}(d\zeta)$ under the map $\zeta \mapsto \cos \zeta$, then  
\begin{equation*}
(j+\nu) \int C_j^{\nu}(\cos \zeta)C_m^{\nu}(\cos \zeta) \mu^{\nu}(d\cos \zeta) = \nu C_j^{\nu}(1) \delta_{jm}
\end{equation*}
so that \eqref{Series} yields
\begin{align}\label{E1}
\nu (\pm 1)^j \left(\frac{2}{R}\right)^{p\nu}I_{p(j+\nu)}(R) = \int P_j^{\nu}(\cos \zeta) f_{\nu,p}^{\pm}(R,\cos \zeta) \mu^{\nu}(d\cos\zeta)
\end{align}
where 
\begin{equation*}
P_j^{\nu}(\cos \zeta) := C_j^{\nu}(\cos \zeta)/C_j^{\nu}(1)
\end{equation*}
is the $j$-th normalized Gegenbauer polynomial. Thus, the $j$-th Gegenbauer-Fourrier coefficients of $\zeta \mapsto f_{\nu,p}^{\pm}(R,\cos \zeta)$ are given by 
\begin{equation*}
\nu (\pm 1)^j \left(\frac{2}{R}\right)^{p\nu} {\it I}_{p(j +\nu)}(R), \quad p \geq 1. 
\end{equation*}
Following \cite{Abou} p.356, the Mehler integral (\cite{Erd}, p.177)
\begin{equation*}
P_j^{\nu}(\cos \zeta) =  2^{\nu}\frac{\Gamma(\nu+1/2)}{\Gamma(\nu)\sqrt{\pi}} (\sin \zeta)^{1-2\nu} \int_0^{\zeta} \cos[(j+\nu)t] (\cos t - \cos \zeta)^{\nu-1} dt
\end{equation*}
valid for real $\nu > 0$, transforms \eqref{E1} to 
\begin{align} 
\left(\frac{2}{R}\right)^{p\nu}(\pm 1)^jI_{p(j+\nu)}(R) &= \frac{2^{\nu}}{\pi} \int_0^{\pi} f_{\nu,p}^{\pm}(R,\cos \zeta)\sin \zeta  \int_0^{\zeta} \cos[(j+\nu)t] (\cos t - \cos \zeta)^{\nu-1} dt d\zeta \nonumber
\\ & = \label{F1} \frac{2^{\nu}}{\pi} \int_0^{\pi} \cos[(j+\nu)t]  \int_{t}^{\pi} f_{\nu,p}^{\pm}(R,\cos \zeta)\sin \zeta  (\cos t - \cos \zeta)^{\nu-1} d\zeta dt.
\end{align}
The second integral displayed in the RHS of the second equality is known as the Radon transform of $\zeta \mapsto f_{\nu,p}^{\pm}(R,\cos \zeta)$ and inversion formulas already exist (\cite{Abou}). As a matter of fact, we firstly need to express $(\pm 1)^{j+\nu} {\it I}_{p(j+\nu)}$, when $\nu \geq 1$ is an integer, as the Fourier-cosine coefficient of order $j+\nu$ of some function. This is a consequence of the Lemma below. Secondly, we shall use the appropriate inversion formula for the Radon transform derived in \cite{Abou}).  

\begin{Lemma}\label{L}
For any integer $p \geq 1$ and any $t \in [0,\pi]$: 
\begin{eqnarray*}
 2 \sum_{j \geq 0} {\it I}_{pj}(R) \cos (jt) & = & {\it I}_0(R) + \frac{1}{p}\sum_{s=1}^{p}e^{R\cos[(t+2\pi s )/p]} \\ 
 2 \sum_{j \geq 0} (-1)^j {\it I}_{pj}(R) \cos (jt) & = & {\it I}_0(R) + \frac{1}{p}\sum_{s=1}^{p}e^{R\cos[(t+2\pi s + \pi)/p]}.
\end{eqnarray*} 
\end{Lemma}

{\it Proof.} 
We will prove the first equality, the proof of the second one follows from $(-1)^j\cos(jt) = \cos(j(t+\pi))$. Write
\begin{align*}
2 \sum_{j \geq 0} {\it I}_{pj}(R) \cos (jt) & = \sum_{j \geq 0} {\it I}_{pj}(R) [e^{ijt} + e^{-ijt}] 
\\& = {\it I}_0(R) + \sum_{j \in \mathbb{Z}} {\it I}_{pj}(R) e^{ijt}
\end{align*} 
where we used the fact that ${\it I}_j(r) = {\it I}_{-j}(r), j \geq 0$. Using the identity
\begin{equation}
\frac{1}{m}\sum_{s=1}^{m} e^{2i\pi sj/m} = 
\left\{\begin{array}{lcr} 
1 & \textrm{if} & j \equiv 0[m], \\
0 & \textrm{otherwise}, & 
\end{array}\right.
\end{equation}
valid for any integer $m \geq 1$, one obviously gets 
\begin{align*}
\sum_{j \in \mathbb{Z}} {\it I}_{pj}(R) e^{ijt} = \frac{1}{p} \sum_{s=1}^p \sum_{j \in \mathbb{Z}} {\it I}_j(R) e^{ij(t+ 2\pi s)/p}. 
\end{align*}
The first equality of the Lemma then follows from the generating series for modified Bessel functions (\cite{Watson}): 
\begin{equation*}
e^{(z + 1/z)R/2} = \sum_{j \in \mathbb{Z}}{\it I}_j(R)z^j, \, z \in \mathbb{C}. \
\end{equation*}
Now, we use the Lemma to get
\begin{equation*}
{\it I}_{pj}(R) = {\it I}_0(R) \delta_{j0} + \frac{1}{\pi}\int_0^{\pi} \cos (jt) \frac{1}{p} \sum_{s=1}^p e^{R\cos[(t+2\pi s)/p]}dt
\end{equation*}
for any integer $j \geq 0$. Assuming that $\nu$ is a strictly positive integer, one then recovers 
\begin{equation}\label{F2}
{\it I}_{p(j+\nu)}(R) =  \frac{1}{\pi}\int_0^{\pi} \cos [(j+\nu)t] \frac{1}{p} \sum_{s=1}^p e^{R\cos[(t+2\pi s)/p]}dt. 
\end{equation}
Note that 
\begin{equation*}
t \mapsto \int_{t}^{\pi} f_{\nu,p}^+(R,\cos \zeta)\sin \zeta  (\cos t - \cos \zeta)^{\nu-1}d\zeta
\end{equation*}
as well as 
\begin{equation*}
t \mapsto \frac{1}{p} \sum_{s=1}^p e^{R\cos[(t+2\pi s)/p]}
\end{equation*}
are even functions. This is true since 
\begin{equation*}
\zeta \mapsto f_{\nu,p}^+(R,\cos \zeta)(\sin \zeta)  (\cos t - \cos \zeta)^{\nu-1}
 \end{equation*}
is an odd function so that 
\begin{equation*}
\int_{-t}^t f_{\nu,p}^+(R,\cos \zeta)\sin \zeta  (\cos t - \cos \zeta)^{\nu-1} d\zeta = 0,
\end{equation*}
and  since 
 \begin{equation*}
 \cos[(-t+2s\pi)/p] = \cos[(t + 2(p-s)\pi)/p]
 \end{equation*}
so that one performs the index change $s \rightarrow p-s$ and notes that the terms corresponding to $s=0$ and $s=p$ are equal. Similar arguments yield the $2\pi$-periodicity of these functions, therefore the Fourier-cosine transforms of their restrictions on $(-\pi,\pi)$ coincide with their Fourier transforms on that interval. As a matter of fact,  
\begin{equation*}
\left(\frac{R}{2}\right)^{p\nu}\int_{t}^{\pi}f_{\nu,p}^+(R,\cos \zeta)\sin \zeta  (\cos t - \cos \zeta)^{\nu-1} d\zeta  =  \frac{1}{2^{\nu} p}\sum_{s=1}^p e^{R\cos[(t+2\pi s)/p]}
\end{equation*}
for all $t$ since  both functions are continuous and similarly
\begin{equation*}
\left(\frac{R}{2}\right)^{p\nu}\int_{t}^{\pi}f_{\nu,p}^-(R,\cos \zeta)\sin \zeta  (\cos t - \cos \zeta)^{\nu-1} d\zeta  =  \frac{(-1)^{\nu}}{2^{\nu} p}\sum_{s=1}^p e^{R\cos[(t+2\pi s+\pi)/p]}
\end{equation*}
Finally, our main result follows from Theorem 3.1. p.363 in \cite{Abou}. 

\subsection{Second proof} 
It uses Lemma \ref{L} together with the formula (\cite{Erd}): 
\begin{equation*}
\frac{1}{2^{\nu-1}\Gamma(\nu)} \frac{d^{\nu}}{dz^{\nu}}T_{j+\nu}(z) = (j+\nu)C_j^{\nu}(z)
 \end{equation*}
valid for integer $\nu \geq 1$, where $T_j$ is the $j$-th Tchbycheff polynomial of the first kind defined by 
\begin{equation*}
T_j(z) = \cos(j\arccos z), \, j \geq 0.
\end{equation*}
Indeed, an elementary change of variables gives
\begin{equation*}
\left[-\frac{1}{\sin \xi} \frac{d}{d\xi}\right]^{\nu}[T_{j+\nu}(\cos (\cdot))](\zeta) = C_j^{\nu}(\cos\zeta), 
\end{equation*}
and a standard argument from analysis allows to write
\begin{align*}
\left(\frac{R}{2}\right)^{p\nu} f_{\nu,p}^+(R,\cos\zeta) &= \frac{1}{2^{\nu-1}\Gamma(\nu)}\left[-\frac{1}{\sin \xi} \frac{d}{d\xi}\right]^{\nu} \sum_{j \geq 0}{\it I}_{p(j +\nu)}(R) T_{j+\nu}(\cos\zeta)
\\& = \frac{1}{2^{\nu-1}\Gamma(\nu)}\left[-\frac{1}{\sin \xi} \frac{d}{d\xi}\right]^{\nu} \sum_{j \geq -\nu}{\it I}_{p(j +\nu)}(R) \cos[(j+\nu)\zeta]
\\& = \frac{1}{2^{\nu-1}\Gamma(\nu)} \left[-\frac{1}{\sin \xi} \frac{d}{d\xi}\right]^{\nu} \sum_{j \geq 0}{\it I}_{pj}(R) \cos(j\zeta)
\\& = \frac{1}{2^{\nu}\Gamma(\nu)} \left[-\frac{1}{\sin \xi} \frac{d}{d\xi}\right]^{\nu} \frac{1}{p}\sum_{s=1}^{p}e^{R\cos[(\zeta+2\pi s)/p]}
\end{align*}
where the second equality follows from the fact that $T_j$ is a polynomial of degree $j$ while the last one follows from Lemma \ref{L}. 


\section{Dihdral groups $D_4,D_6$}
\subsection{p=2}
Letting $p=2$ and using the fact that $\cosh$ is an even function, our main result yields 
\begin{align*}
\left(\frac{4}{R^2}\right)^{\nu}\sum_{j \geq 0}(2j+\nu){\it I}_{2(2j +\nu)}(R)C_{2j}^{\nu}(\cos \zeta) = \frac{1}{2^{\nu}\Gamma(\nu)} \left[\frac{4}{R^2\sin \zeta} \frac{d}{d\zeta}\right]^{\nu} \\ 
\left[(-1)^{\nu}\cosh \left(R\cos (\cdot/2)\right) + \left(\cosh(R\sin(\cdot/2)\right)\right] (\zeta).
\end{align*}
Noting that for a function $f$
\begin{equation*}
-\frac{4}{R^2\sin \zeta} \frac{d}{d\zeta} f\left(R\cos \frac{\zeta}{2}\right) (\zeta) = \left[\frac{1}{u} \frac{d}{du}f(u)\right]_{|u = R\cos (\zeta/2)}, 
\end{equation*}
\begin{equation*}
\frac{4}{R^2\sin \zeta} \frac{d}{d\zeta} f\left(R\sin \frac{\zeta}{2}\right) (\zeta) = \left[\frac{1}{u} \frac{d}{du}f(u)\right]_{|u = R\sin (\zeta/2)}, 
\end{equation*}
and using the classical formula (see for instance (5.8.3) in \cite{Leb})
\begin{align*}
\left(\frac{1}{z}\frac{d}{dz}\right)^{\nu-1} \frac{\sin z}{z} = (-1)^{\nu-1}\sqrt{\frac{\pi}{2}} \frac{1}{z^{\nu-1/2}} J_{\nu-1/2}(z),
\end{align*}
one obtains  
\begin{equation*}
\left(\frac{4}{R^2}\right)^{\nu}\sum_{j \geq 0}(2j+\nu){\it I}_{2(2j +\nu)}(R)C_{2j}^{\nu}(\cos \zeta) = \frac{1}{2\Gamma(2\nu)}  \left[i_{\nu -1/2}\left(R\cos \frac{\zeta}{2}\right) +  i_{\nu -1/2}\left(R\sin \frac{\zeta}{2}\right)\right],
\end{equation*}
Finally, one recovers Corollary 1.2 in \cite{Dem02} since 
\begin{equation*}
\frac{c_{2,k}}{c(k_1-1/2,k_0-1/2)} = \frac{\Gamma(2\nu+1)}{\nu},
\end{equation*}
and using the known formulas
\begin{equation*}
\cos(\zeta/2) = \sqrt{\frac{1+\cos \zeta}{2}}, \quad \sin(\zeta/2) = \sqrt{\frac{1-\cos \zeta}{2}}.
\end{equation*}

\subsection{p=3}
The corresponding dihedral group ${\it D}_6$ is a two dimensional representation of the exceptional Weyl group $G_2$ (\cite{Baez}). Let $\zeta \in ]0,\pi[$ and start with the linearization formula: 
\begin{equation*}
4 \cos^3(\zeta/3) = \cos \zeta + 3\cos (\zeta/3).
\end{equation*}
Thus, we are led to find a root lying in $[-1,1]$ of the cubic equation 
\begin{equation*}
Z^3 - (3/4)Z - (\cos \zeta)/4 = 0 
\end{equation*}
for $|Z| < 1$. Set $Z = (\sqrt{-1}/2)T, |T| < 2$, the above cubic equation transforms to 
\begin{equation*}
T^3 + 3T -2\sqrt{-1}\cos \zeta = 0.  
\end{equation*}
The obtained cubic equation already showed up in analytic function theory in relation to the local inversion Theorem (\cite{Hille} p.265-266). Amazingly (compared to Cardan formulas), its real and both complex roots are expressed through the Gauss hypergeometric function ${}_2F_1$. Since we are looking for real $Z = (\sqrt{-1}/2)T$, we shall only consider the complex roots (see the bottom of p. 266 in \cite{Hille}): 
\begin{equation*}
T^{\pm} = \pm \sqrt{-1} \left[\sqrt{3}\, {}_2F_1\left(-\frac{1}{6},\frac{1}{6}, \frac{1}{2}; \cos^2 \zeta\right) - \frac{1}{3}\cos \zeta \, {}_2F_1\left(\frac{1}{3},\frac{2}{3}, \frac{3}{2}; \cos^2 \zeta\right)
\right]
\end{equation*}
so that
\begin{equation*}
Z^{\pm}  = \pm \left[\frac{\sqrt{3}}{2}\, {}_2F_1\left(-\frac{1}{6},\frac{1}{6}, \frac{1}{2}; \cos^2 \zeta\right) - \frac{1}{6}\cos \zeta \, {}_2F_1\left(\frac{1}{3},\frac{2}{3}, \frac{3}{2}; \cos^2 \zeta\right)
\right]. 
\end{equation*}
Since for $\zeta = \pi/2$, $\cos \zeta/3 = \cos \pi/6 = \sqrt{3}/2$, it follows that
\begin{equation*}
\cos(\zeta/3) = \left[\frac{\sqrt{3}}{2}\, {}_2F_1\left(-\frac{1}{6},\frac{1}{6}, \frac{1}{2}; \cos^2 \zeta\right) - \frac{1}{6}\cos \zeta \, {}_2F_1\left(\frac{1}{3},\frac{2}{3}, \frac{3}{2}; \cos^2 \zeta\right)
\right]
\end{equation*}
for all $\zeta \in (0,\pi)$. Now, write $Z = Z(\cos\zeta)$ so that 
\begin{align*}
\cos[(\zeta+2s\pi)/3] &= \cos(2s\pi/3)\cos(\zeta/3) - \sin(2s\pi/3)\sqrt{1-\cos^2(\zeta/3)}\\&
=  \cos(2s\pi/3)Z(\cos \zeta) - \sin(2s\pi/3)\sqrt{1-Z^2(\cos\zeta)}
\end{align*}
for any $1 \leq s \leq 3$. It follows that
\begin{equation*}
f_{\nu,3}^+(R,\cos\zeta) = \frac{1}{3\Gamma(\nu)} \left[-\frac{4}{R^3\sin \zeta} \frac{d}{d\zeta}\right]^{\nu}  \sum_{s=1}^3g_s(RZ(\cos\zeta))
 \end{equation*}
where 
\begin{equation*}
g_s(u) = \exp\{\left[\left(\cos(2s\pi/3)u - \sin(2s\pi/3)\sqrt{R^2-u^2}\right)\right]\}, u \in (-1,1).
\end{equation*}
Equivalently 
\begin{equation*}
f_{\nu,3}^+(R,\cos\zeta) = \frac{1}{3\Gamma(\nu)} \left[\frac{4}{R^3} \frac{d}{du}\right]^{\nu}\sum_{s=1}^3h_s(u)_{|u=\cos \zeta}
\end{equation*}
where $h_s(u) := g_s(RZ(u)), 1 \leq s \leq 3$. For instance, let $\nu = 1$, then it is not difficult to see that 
\begin{equation*}
\frac{d}{du} h_s(u)_{|u=\cos \zeta} = \frac{R}{\sin \zeta/3}\frac{dZ}{du}_{|u=\cos \zeta} \sin\left(\frac{\xi + 2\pi s}{3}\right) \sinh\left[\sin\left(\frac{\xi + 2\pi s}{3}\right)\right]
\end{equation*}
for any $s\in \{1,2,3\}$ and the derivative of $u \mapsto Z(u)$ is computed using the differentiation formula for ${}_2F_1$: 
\begin{equation*}
\frac{d}{du} {}_2F_1(a,b,c;u) = \frac{ab}{c} {}_2F_1(a+1,b+1,c+1;u), \, |u| < 1, c \neq 0. 
\end{equation*}
Similar results may be derived for $f_{\nu,3}^-$ and the reader may realize that formulas are cumbersome compared to the ones derived for $p=2$. Nonetheless, computations seems to be still tractable for $p=4$. 


\subsection{General values of $p$} Let $\xi \in [0,\pi]$ and $ p\geq 1$ then $T_p$ is invertible on the interval $[\cos(\pi/p),1]$. This is easily seen from the very definition of $T_p$: 
\begin{equation*}
\frac{dT_p}{dz}(z)= \frac{1}{\sqrt{1-z^2}}\sin(\arccos(z)) > 0.
\end{equation*}  
Consequently, since $T_p[\cos(\zeta/p)] = \cos\xi$ then 
\begin{align*}
\left(\frac{R}{2}\right)^{p\nu} f_{\nu,p}^+(R,\cos\zeta) = \frac{1}{2^{\nu}\Gamma(\nu)} \left[-\frac{1}{\sin \zeta} \frac{d}{d\zeta}\right]^{\nu} \frac{1}{p}\sum_{s=1}^{p}e^{R[\cos(2\pi s/p) T_p^{-1}(\cos\zeta) - \sin(2\pi s/p) \sqrt{1-T_p^{-2}(\cos\zeta)}]}
\end{align*}
where we used the notation $T_p^{-2} = (T_p^{-1})^2$. It follows that 
\begin{align*}
\left(\frac{R}{2}\right)^{p\nu} f_{\nu,p}^+(R,\cos \zeta) = \frac{1}{2^{\nu}\Gamma(\nu)} \frac{1}{p}\sum_{s=1}^{p} \left[\frac{d}{dz}\right]^{\nu} e^{R[\cos(2\pi s/p) T_p^{-1}(z) - \sin(2\pi s/p) \sqrt{1-T_p^{-2}(z)}]}{|}_{z=\cos \zeta},
\end{align*}
and similarly 
\begin{align*}
\left(\frac{R}{2}\right)^{p\nu} f_{\nu,p}^-(R,\cos \zeta) = \frac{1}{2^{\nu}\Gamma(\nu)} \frac{1}{p}\sum_{s=1}^{p} \left[\frac{d}{dz}\right]^{\nu} e^{R[\cos(2\pi s/p) T_p^{-1}(z) - \sin(2\pi s/p) \sqrt{1-T_p^{-2}(z)}]}{|}_{z=-\cos \zeta}.
\end{align*}
However both $f_{\nu,p}^{\pm}$ give the same contribution in the expression of $D_k^W$ due to the integration with respect to symmetric Beta distributions and since 
\begin{equation*}
\pm \cos\xi = \pm z_{p\phi,p\theta}(u,v) = \pm[u\cos(p\phi) \cos(p\theta) + v \sin(p\phi) \sin(p\theta)]. 
\end{equation*}

\section{Odd Dihedral groups}
Let $n \geq 3$ be an odd integer and consider odd dihedral groups ${\it D}_n$, then (\cite{Dem01} p.157) 
\begin{equation}\label{GBFE1}
D_k^W(\rho,\phi,r,\theta) = c_{n,k}\left(\frac{2}{r\rho}\right)^{nk} \sum_{j \geq 0}{\it I}_{n(2j+k)}(\rho r)p_j^{-1/2, l_0}(\cos(2n\phi))p_j^{-1/2, l_0}(\cos(2n\theta))
\end{equation}
where $k \geq 0, \rho,r \geq 0, \theta,\phi \in [0,\pi/n]$, and 
\begin{equation*}
c_{n,k} = 2^{k}\Gamma(nk+1)\frac{\sqrt{\pi}\Gamma(k+1/2)}{\Gamma(k+1)}.
\end{equation*}
In order to adapt our method to these groups, we need to write down the product formula for orthonormal Jacobi polynomials in the limiting case $\alpha = -1/2$ or equivalently $k=0$. This task was achieved in \cite{Dij} p.194 using implicitly the fact that the Beta distribution $\mu^{\alpha}$ converges weakly to the Dirac mass $\delta_1$ as $\alpha \rightarrow -1/2$. In order to fit it into our normalizations, we proceed as follows: use 
 the well-known quadratic transformation (\cite{Dunkl}): 
\begin{align*}
P_j^{-1/2,k-1/2}(1-2\sin^2(n\theta)) &= (-1)^j P_j^{k-1/2,-1/2}(2\sin^2(n\theta) - 1) 
\\& = (-1)^j \frac{(1/2)_j}{(k)_j} C_{2j}^k(\sin(n\theta))
\end{align*}
where $P_j^{\alpha,\beta}$ is the (non orthonormal) $j$-th Jacobi polynomial, together with $\cos(2n\theta) = 1-2\sin^2(n\theta)$ to obtain 
\begin{equation*}
P_j^{-1/2,k-1/2}(\cos(2n\theta))P_j^{-1/2,k-1/2}(\cos(2n\phi)) = \left[\frac{(1/2)_j}{(k)_j}\right]^2 C_{2j}^k(\sin(n\theta))C_{2j}^k(\sin(n\phi)).
\end{equation*}
Now, let $k > 0$ and recall that the squared $L^2$-norm of $P_j^{-1/2,k-1/2}$ is given by (\cite{Dunkl})
\begin{align*}
\frac{2^k}{2j+k} \frac{\Gamma(j+1/2)\Gamma(j+k+1/2)}{j!\Gamma(j+k)} = \frac{2^k\sqrt{\pi}\Gamma(k+1/2)}{\Gamma(k)} \frac{(1/2)_j}{(k)_j} \frac{(k+1/2)_j}{(2j+k)j!}.
\end{align*}
Recall also the special value 
\begin{align*}
C_{2j}^k(1) = \frac{(2k)_{2j}}{(2j)!} 
=2\frac{(k)_j(k+1/2)_j}{(1/2)_jj!}.
\end{align*}
It follows that 
\begin{align*}
c(k)p_j^{-1/2,k-1/2}(\cos(2n\theta))p_j^{-1/2,k-1/2}(\cos(2n\phi)) &=  \frac{(1/2)_j}{(k)_j}\frac{(2j+k)j!}{(k+1/2)_j}C_{2j}^k(\sin(n\theta))C_{2j}^k(\sin(n\phi))
\\&=  \frac{(2j+k)}{C_{2j}^k(1)}C_{2j}^k(\sin(n\theta))C_{2j}^k(\sin(n\phi))
\\& = (2j+k) \int C_{2j}^k\left(z_{n\phi,n\theta}(u,1)\right)\mu^k(du),
\end{align*}
according to \cite{Dij} p.194, where 
\begin{equation*}
c(k) := \frac{2^{k+1}\sqrt{\pi}\Gamma(k+1/2)}{\Gamma(k)}.
\end{equation*}
As a matter of fact, we are led again to 
\begin{equation*}
\left(\frac{2}{R}\right)^{nk} \sum_{j \geq 0}(2j+k){\it I}_{n(2j +k)}(R)C_{2j}^{k}(\cos \zeta) = \frac{1}{2}[f_{k,n}^++ f_{k,n}^-](R,\cos \zeta). 
\end{equation*}

\section{Concluding Remarks}
1/The two proofs we wrote in this paper come somehow in opposite ways. More precisely, the operator 
\begin{equation*}
\left[-\frac{1}{\sin \zeta}\frac{d}{d\zeta}\right]^{\nu} 
\end{equation*}
appears at the second step of the first proof, that is when inverting the Radon transform of the function 
\begin{equation*}
\sum_{j \geq 0}{\it I}_{pj}(R)\cos(j\zeta),
\end{equation*}
while it appears at the beginning of the second proof due to the existing relation between $T_{j+\nu}$ and $C_j^{\nu}$ for integer $\nu$.  \\
2/Our main result explains the increasing disability of getting a `nice' expression of the generalized Bessel function associated with dihedral groups as the integer parameter $p$ increases. This disability comes in inverting the $p$-th Tchebycheff polynomial of the first kind in the interval $[\cos(\pi/p), 1]$, therefore finding a suitable root of a polynomial of degree $p$. \\
3/In \cite{Abou}, authors defined and studied the so-called Radon transform on an Euclidean sphere $S^d \in \mathbb{R}^{d+1}, d \geq 2$. It is given, up to a factor, by an average over circles passing through a fixed point on a this sphere  (see \cite{Abou} (15), p. 359 and (29), p. 369) and we easily identify $\nu = (d-1)/2$ (see (10), p. 357). Consequently $d = 2k+1$ for odd dihedral groups, a fact that is related to $D_k^W$ as follows. Recall that $D_k^W(\cdot,y)$ is an eigenfunction of  the Dunkl-Laplace operator and that the latter acts on ${\it D}_n$-invariant functions as (\cite{Dem01}): 
\begin{equation*}
\Delta_k^W  :=  \partial_r^2  + \frac{2nk+1}{r}\partial_r  + \frac{1}{r^2}\left[\partial_{\theta}^2 + 2nk\cot(n\theta)\partial_{\theta}\right]
\end{equation*}
for odd integers $n$. Now transform the `angular' part of $\Delta_k^W$ 
\begin{equation*}
\partial_{\theta}^2 + 2nk\cot(n\theta)\partial_{\theta} 
\end{equation*}
to 
\begin{equation*}
\partial_{\theta}^2 + 2k\cot(\theta)\partial_{\theta}
\end{equation*}
to see that it fits the Euclidean Laplacian acting on $SO(2k+1)$-invariant functions on the unit sphere $S^{2k+1} \approx SO(2k+2)/SO(2k+1)$ (see Proposition 2.3 p.197 in \cite{Faraut}). But the spherical functions of the Gelfand pair $(SO(2k+2), SO(2k+1))$ are expressed through Gegenbauer polynomials of index $[(2k+2)-2]/2 = k = \nu$ (see \cite {Abou} p. 356). A similar statement holds for even dihedral groups and equal multiplicities $k_0 = k_1$: 
the Dunkl-Laplace operator acts on ${\it D}_{2p}$-invariant functions as
\begin{equation*}
\Delta_k^W = \partial_r^2  + \frac{2p(k_0+k_1) +1}{r}\partial_r + \frac{1}{r^2}\left[\partial_{\theta}^2+ 2p(k_0\cot(p\theta) - k_1\tan(p\theta))\partial_{\theta}\right]
\end{equation*}
and the following formula is obvious
\begin{equation*}
\cot(\theta) - \tan(\theta) = 2\cot(2\theta). 
\end{equation*} 
4/We learnt from our colleague S. Ben Said that $f_{\nu,1}^+$ and $f_{\nu, 2}^+$ already appeared in relation to representation theory of the metapletic and the indefinite orthogonal groups (\cite{KobMan}). We do not know whether $f_{\nu,p}^+$ is related for general $p$ to the representation theory of some groups.  
 
{\bf Acknowledgment}: the author is grateful to Professor C.F. Dunkl who made him aware of the hypergeometric formulas for the roots of the cubic equation.       


\end{document}